\documentclass[11pt]{amsart}
\usepackage{mathrsfs}
\usepackage{amssymb,latexsym}
\usepackage{pstcol,pstricks,color}

 \numberwithin{equation}{section}

\newtheorem{theorem}{Theorem}[section]

\newtheorem{corollary}[theorem]{Corollary}

\newtheorem{conjecture}[theorem]{Conjecture}
\newtheorem{remark}[theorem]{Remark}

\newtheorem{lemma}[theorem]{Lemma}

\def\qed{\hfill $\Box$}
\def\pf{\noindent {\it Proof.} }

\title{The largest singletons of set partitions}

\begin{document}
\maketitle
\begin{center}
Yidong Sun$^\dag$ and Xiaojuan Wu


Department of Mathematics, Dalian Maritime University, 116026 Dalian, P.R. China\\[5pt]

{\it $^\dag$Email: sydmath@yahoo.com.cn   }

{\bf Dedicated to L.C. Hsu, on the occasion of his ninetieth birthday}
\end{center}\vskip0.2cm

\subsection*{Abstract}
Recently, Deutsch and Elizalde studied the largest and the smallest fixed points of
permutations. Motivated by their work, we consider the analogous problems in set partitions.
Let $A_{n,k}$ denote the number of partitions of $\{1,2,\dots, n+1\}$ with the largest singleton $\{k+1\}$ for $0\leq k\leq n$.
In this paper, several explicit formulas for $A_{n,k}$, involving
a Dobinski-type analog, are obtained by algebraic and combinatorial methods,
many combinatorial identities involving $A_{n,k}$ and Bell numbers are presented by
operator methods, and congruence properties of $A_{n,k}$ are also investigated. It will been showed that the sequences
$(A_{n+k,k})_{n\geq 0}$ and $(A_{n+k,k})_{k\geq 0}$ (mod $p$) are periodic for any prime $p$, and contain a string of $p-1$ consecutive zeroes.
Moreover their minimum periods are conjectured to be
$N_p=\frac{p^p-1}{p-1}$ for any prime $p$.

\medskip

{\bf Keywords}: Set partition; Bell number; Congruence; Identity.

\noindent {\sc 2000 Mathematics Subject Classification}:
05A05; 05A18; 05A19; 05A40; 11A07.

\section{Introduction}

A {\it partition} of a set $[n]=\{1, 2, \dots, n\}$ is a
collection of nonempty and mutually disjoint subsets of $[n]$,
called {\it blocks}, whose union is $[n]$. It is well known that the number of partitions
of $[n]$ with exactly $k$ blocks is the Stirling number of the
second kind $S(n, k)$ \cite[A008267]{Sloane} and the total number of partitions of $[n]$ is
the $n$-th Bell number $B_n$ \cite{Rota}, beginning with $(B_n)_{n\geq 0}=(1,1,2,5,15,52,203,\dots)$ \cite[A000110]{Sloane} and having
the exponential generating function \cite{Stanley}
\begin{eqnarray}\label{eqn 1.1}
B(x)=\sum_{n\geq 0}B_{n}\frac{x^n}{n!}=\exp(e^{x}-1).
\end{eqnarray}
Differentiating (\ref{eqn 1.1}) gives $B'(x) = e^{x}B(x)$, which leads to
\begin{eqnarray}\label{eqn 1.2}
B_{n+1}=\sum_{k=0}^{n}\binom{n}{k}B_{k}.
\end{eqnarray}

A {\it singleton} of a partition is a block containing just one element. If $\{k\}$ is a singleton of a partition, we
denote it by $k$ for short. The number of partitions
of $[n]$ without singletons is counted by $V_n$ beginning with $(V_n)_{n\geq 0}=(1,0,1,1,4,11,41,162,\dots)$
\cite[A000296]{Sloane}, and having the exponential generating function
\begin{eqnarray}\label{eqn 1.4}
V(x)=\sum_{n\geq 0}V_{n}\frac{x^n}{n!}=\exp(e^{x}-x-1).
\end{eqnarray}
Bernhart \cite{Bernhart} has given a combinatorial interpretation for the relation $B_n=V_n+V_{n+1}$ which can also be obtain from
$B(x)=V(x)+V'(x)$. By (\ref{eqn 1.1}) and (\ref{eqn 1.4}), one can deduce that
\begin{eqnarray*}
B_{n}=\sum_{j=0}^{n}\binom{n}{j}V_j \hskip.5cm {\mbox{and}} \hskip.5cm V_{n}=\sum_{j=0}^{n}(-1)^{n-j}\binom{n}{j}B_j.
\end{eqnarray*}

Recently, Deutsch and Elizalde \cite{DeuEliz} study the largest and the smallest fixed points of
permutations. Motivated by their work, we consider the
analogous problems in set partitions. Let $A_{n,k}$ denote the number of partitions of $[n+1]$ with the largest singleton $k+1$.
Clearly,
\begin{eqnarray*}
A_{n,0}=V_{n} \hskip.5cm {\mbox{and}} \hskip.5cm A_{n,n}=B_{n}.
\end{eqnarray*}

This paper is organized as follows. In the next section, we find several explicit formulas for $A_{n,k}$,
involving a Dobinski-type analog, by algebraic and combinatorial methods. In the section 3, we obtain
many combinatorial identities involving $A_{n,k}$ and Bell numbers $B_n$ by operator methods. In the last section,
we consider the congruence properties of $A_{n,k}$ and Bell numbers $B_n$, find that the sequences
$(A_{n+k,k})_{n\geq 0}$ and $(A_{n+k,k})_{k\geq 0}$ (modulo $p$) are periodic for any prime $p$ and contain a string of $p-1$ consecutive zeroes.
We also conjecture that their minimum periods are $N_p=\frac{p^p-1}{p-1}$ for any prime $p$.

\section{The explicit formulas for $A_{n,k}$}

It follows from the definition that
\begin{eqnarray}\label{eqn 1.5}
A_{n,k}=V_n+\sum_{j=0}^{k-1}A_{n-1,j},
\end{eqnarray}
since by removing the largest singleton $k+1$ of a partition of $[n+1]$ containing singletons, we get a partition of
$\{1, \dots, k, k+2, \dots, n+1\}$ whose largest singleton (if any) is less than $k+1$.

In (\ref{eqn 1.5}), if we replace $k$ by $k-1$, then by subtraction we obtain
a recurrence for $n, k\geq 1$,
\begin{eqnarray}\label{eqn 2.1}
A_{n,k}=A_{n,k-1}+A_{n-1, k-1}.
\end{eqnarray}

Table 1 shows the values of $A_{n,k}$ for small $n$ and $k$. It should be noticed that $\{A_{n+k,k}\}_{n\geq k\geq 1}$ is just the Aitken's array
\cite[A011971]{Sloane}. We point out that it is possible to give a direct combinatorial
proof of the recurrence (\ref{eqn 2.1}) from the definition of the $A_{n,k}$. Indeed, given a partition $\pi$ of
$[n+1]$ with the largest singleton $k+1$, if $k$ is also a singleton, delete the singleton $k+1$ and subtracting one from all the entries
large than $k+1$, we obtain a partition of $[n]$ with the largest singleton $k$; if $k$ is not a singleton, exchange $k$ and $k+1$,
we obtain a partition of $[n+1]$ with the largest singleton $k$.

\begin{center}
\begin{eqnarray*}
\begin{array}{c|cccccccc}\hline
n/k & 0   & 1   & 2    & 3    & 4    & 5    & 6     & 7 \\\hline
  0 & 1   &     &      &      &      &      &       &    \\
  1 & 0   & 1   &      &      &      &      &       &   \\
  2 & 1   & 1   & 2    &      &      &      &       &   \\
  3 & 1   & 2   & 3    & 5    &      &      &       &   \\
  4 & 4   & 5   & 7    & 10   &  15  &      &       &    \\
  5 & 11  & 15  & 20   & 27   &  37  &  52  &       &    \\
  6 & 41  & 52  & 67   & 87   &  114 &  151 &  203  &     \\
  7 & 162 & 203 & 255  & 322  &  409 &  523 &  674  & 877 \\\hline
\end{array}
\end{eqnarray*}
Table 1. The values of $A_{n,k}$ for $n$ and $k$ up to $7$.
\end{center}
When $k=1$, (\ref{eqn 2.1}) produces a new setting for Bell numbers, namely $A_{n+1,1}=B_n$. A simple combinatorial proof reads:
given a partition $\pi$ of $[n+2]$ with the largest singleton $2$, if $1$ is also a singleton, delete the two singletons $1, 2$
and subtracting two from all the entries large than $2$, we obtain a partition of $[n]$ without singletons;
if $1$ is not a singleton, break the block containing $1$ into singletons (more than one), then delete the two singletons $1, 2$
and subtracting two from all the entries large than $2$, we obtain a partition of $[n]$ with singletons.

\begin{lemma}\label{lemma 2.1}
The bivariate exponential generating function for $A_{n+k,k}$ is
given by
\begin{eqnarray*}
A(x, y)=\sum_{n,k\geq
0}A_{n+k,k}\frac{x^n}{n!}\frac{y^k}{k!}=\exp(e^{x+y}-x-1).
\end{eqnarray*}
\end{lemma}
\pf Define
\begin{eqnarray*}
A_k(x)=\sum_{n\geq
0}A_{n+k,k}\frac{x^n}{n!}.
\end{eqnarray*}
Clearly, $A_0(x)=\exp(e^x-x-1)$ and $A_1(x)=\exp(e^x-1)$. From (\ref{eqn 2.1}), one can derive that
\begin{eqnarray*}
A_k(x)=A_{k-1}(x)+A_{k-1}'(x).
\end{eqnarray*}
Let $\mathcal{D}$ denote the derivative with respect to $x$, we have
\begin{eqnarray*}
A_k(x)=(1+\mathcal{D})A_{k-1}(x)=(1+\mathcal{D})^kA_0(x).
\end{eqnarray*}
Then
\begin{eqnarray*}
A(x,y) &=& \sum_{k\geq 0}A_k(x)\frac{y^k}{k!}=\sum_{k\geq 0}\frac{y^k(1+\mathcal{D})^k}{k!}A_0(x) \\
       &=& e^{y+y\mathcal{D}}A_0(x)=e^{y}e^{y\mathcal{D}}A_0(x)=e^{y}A_0(x+y) \\
       &=& \exp(e^{x+y}-x-1).
\end{eqnarray*}
This complete the proof. \qed\vskip.2cm

The general formula for the Bell polynomial $B_k(x)=\sum_{j=0}^{k}S(k,j)x^{j}$ states that
\begin{eqnarray*}
B_k(x)=e^{-x}\sum_{m\geq 0}\frac{m^{k}x^m}{m!},
\end{eqnarray*}
which, when $x=1$, produces the Dobinski's formula \cite{Rota} for Bell numbers
\begin{eqnarray*}
B_k=\frac{1}{e}\sum_{m\geq 0}\frac{m^{k}}{m!}.
\end{eqnarray*}

Analogously, we can derive a Dobinski-type formula for $A_{n+k,k}$.
\begin{theorem} For any integers $n, k\geq 0$, there holds
\begin{eqnarray}\label{eqn 2.2}
A_{n+k,k}=\frac{1}{e}\sum_{m=0}^{\infty}\frac{m^k(m-1)^{n}}{m!}.
\end{eqnarray}
\end{theorem}
\pf By Lemma \ref{lemma 2.1}, one has
\begin{eqnarray*}
A(x,y) &=& \exp(e^{x+y}-x-1) \\
       &=& e^{-x-1}\sum_{m\geq 0}\frac{e^{(x+y)m}}{m!} \\
       &=& e^{-1}\sum_{m\geq 0}\frac{1}{m!}\sum_{n\geq 0}\frac{(m-1)^nx^n}{n!}\sum_{k\geq 0}\frac{m^ky^k}{k!}  \\
       &=& e^{-1}\sum_{n,k\geq 0}\frac{x^n}{n!}\frac{y^k}{k!}\sum_{m\geq 0}\frac{m^k(m-1)^{n}}{m!},
\end{eqnarray*}
which leads to (\ref{eqn 2.2}) by comparing the coefficients of $\frac{x^n}{n!}\frac{y^k}{k!}$. \qed\vskip.2cm

\begin{remark} According to the Dobinski-type formula for $A_{n+k,k}$, one can deduce the column generating function
$A_k(x)=V(x)B_k(e^x).$ By attracting the coefficient of $\frac{y^k}{k!}$ from $A(x, y)$, one can also find
$A_k(x)=e^{-x}\sum_{n\geq 0}B_{n+k}\frac{x^n}{n!}=e^{-x}\mathcal{D}^{k}B(x).$ Then one has the relation for Bell polynomials
$\mathcal{D}^{k}B(x)=B(x)B_k(e^x)$.
\end{remark}

\begin{theorem} For any integers $n,m, k\geq 0$, there hold
\begin{eqnarray}
A_{n+m,m}&=& \sum_{j=0}^{n}(-1)^{n-j}\binom{n}{j}B_{m+j}, \label{eqn 2.4}\\
A_{n+m+k,m+k}&=& \sum_{j=0}^{m}\binom{m}{j}A_{n+k+j,k}. \label{eqn 2.5}
\end{eqnarray}
\end{theorem}
\pf Note that $A(x, y)=B(x+y)e^{-x}$ and $\frac{\partial^k}{\partial y^k} A(x,y)=A_k(x+y)e^{y}$ from Lemma \ref{lemma 2.1},
by equating the coefficients of
$\frac{x^ny^m}{n!m!}$ in the resulting series,
one can easily deduce (\ref{eqn 2.4})-(\ref{eqn 2.5}). Here we provide a combinatorial proof.

(1) Let $\mathbb{S}$ denote the set of partitions of $[n+m+1]$ containing at least the singleton $m+1$,
Clearly, $|\mathbb{S}|=B_{m+n}$. Let $\mathbb{S}_i$ be
the subset of $\mathbb{S}$ containing another singleton $m+i+1$ for $1\leq i\leq n$. Set $\overline{\mathbb{S}}_i=\mathbb{S}-\mathbb{S}_i$,
then $\bigcap_{i=1}^n\overline{\mathbb{S}}_i$, counted by $A_{n+m,m}$, is just the set of partitions of $[n+m+1]$ with the largest singleton $m+1$.
For any nonempty $(n-j)$-subset $\mathbb{A}\in [n]$, $\bigcap_{i\in \mathbb{A}}\mathbb{S}_i$, counted by $B_{m+j}$,
is the set of partitions of $[n+m+1]$ containing at least the number $n-j+1$ of singletons $m+1$ and $m+i+1$ for all $i\in \mathbb{A}$.
By the Inclusion-Exclusion principle, we have
\begin{eqnarray*}
|\bigcap_{i=1}^n\overline{\mathbb{S}}_i| &=& |\mathbb{S}-\bigcup_{i=1}^{n}\mathbb{S}_i| \\
                                         &=& |\mathbb{S}|+\sum_{j=0}^{n-1}(-1)^{n-j}\binom{n}{j}|\bigcap_{i\in \mathbb{A}, |\mathbb{A}|=n-j}\mathbb{S}_i| \\
                                         &=& \sum_{j=0}^{n}(-1)^{n-j}\binom{n}{j}B_{m+j},
\end{eqnarray*}
which proves (\ref{eqn 2.4}).

(2) A partition $\pi$ of $[n+m+k+1]$ with the largest singleton $m+k+1$ can be obtained as follows.
Suppose that $\pi$ has exactly $m-j$ singletons in $\{k+1, \dots, k+m\}$, there are $\binom{m}{j}$ ways to do this, so the remainder $j$
elements in $\{k+1, \dots, k+m\}$ can not be singletons in $\pi$. These $j$ elements can be regarded as the roles that greater than $m+k+1$, there
are $A_{n+k+j,k}$ ways to produce a partition $\pi'$ of the remainder $n+k+j+1$ elements with the largest singleton $m+k+1$,
then $\pi'$ together with the $m-j$ singletons forms the desired partition $\pi$. Thus there are $\binom{m}{j}A_{n+k+j,k}$ of such
partitions. Summing up all the possible cases yields (\ref{eqn 2.5}).  \qed\vskip.2cm

The cases $k=0$ and $k=1$ in (\ref{eqn 2.5}) produce
\begin{corollary} For any integers $n,m \geq 0$, there hold
\begin{eqnarray}
A_{n+m,m}&=& \sum_{j=0}^{m}\binom{m}{j}V_{n+j}, \nonumber  \\
A_{n+m+1,m+1}&=& \sum_{j=0}^{m}\binom{m}{j}B_{n+j}. \label{eqn 2.3}
\end{eqnarray}
\end{corollary}

\begin{remark}
The case $m:=m+1$ in (\ref{eqn 2.4}), together with (\ref{eqn 2.3}), produces another identity for Bell numbers
\begin{eqnarray*}
\sum_{j=0}^{n}(-1)^{n-j}\binom{n}{j}B_{m+j+1} &=& \sum_{j=0}^{m}\binom{m}{j}B_{n+j}.
\end{eqnarray*}
\end{remark}

Spivey \cite{Spivey} finds a generalized recurrence for Bell numbers
\begin{eqnarray*}
B_{n+k}=\sum_{r=0}^{n}\sum_{j=0}^{k}\binom{n}{r}B_rS(k,j)j^{n-r},
\end{eqnarray*}
and gives it a simple combinatorial proof. This recurrence has been generalized by Belbachir and
Mihoubi \cite{BelMih}, Gould and Quaintance \cite{GouldQuain}.
We also have a similar formula for $A_{n+k,k}$.

\begin{theorem} For any integers $n, k\geq 0$, there hold
\begin{eqnarray}
A_{n+k,k}&=& \sum_{r=0}^{n}\sum_{j=0}^{k}\binom{n}{r}V_{r}S(k,j)j^{n-r}, \label{eqn 2.6}\\
A_{n+k,k}&=& \sum_{r=0}^{n}\sum_{j=0}^{k}\binom{n}{r}B_{r}S(k,j)(j-1)^{n-r}. \label{eqn 2.7}
\end{eqnarray}
\end{theorem}
\pf Note that $A_k(x)=V(x)B_k(e^x)$ and $A_k(x)=B(x)B_k(e^x)e^{-x}$ from Remark 2.3,
by equating the coefficients of
$\frac{x^n}{n!}$ in the resulting series,
one can easily deduce (\ref{eqn 2.6})-(\ref{eqn 2.7}). Here we provide a combinatorial proof.

For the set $[n+k+1]$, one can count the number of ways
to partition these $n+k+1$ elements in the following manners.

(1) Partition the set $[k]$ into
exactly $j$ blocks, there are $S(k,j)$ ways to do this. Choose an $r$-subset from the set $\{k+2, \dots, n+k+1\}$
to be partitioned into new blocks, and distribute the remainder $n-r$ elements among
the $j$ blocks formed from the set $[k]$. There are $\binom{n}{r}$ ways to choose the $r$ elements, $V_r$
ways to partition them into new blocks without singletons, and $j^{n-r}$ ways to distribute the remainder $n-r$
elements among the $j$ blocks. Thus there are $\binom{n}{r}V_{r}S(k,j)j^{n-r}$ of such partitions. Note that $k+1$ is always
a singleton, summing over all possible values of $j$ and $r$ produces all ways to partition the set $[n+k+1]$ with the
largest singleton $k+1$. This gives a proof of (\ref{eqn 2.6}).

(2) Partition the set $[k]$ into
exactly $j$ blocks $\mathbb{S}_1, \mathbb{S}_2,\dots,\mathbb{S}_j$ and assume that $1\in \mathbb{S}_1$, there are $S(k,j)$ ways to do this.
Choose an $r$-subset $\mathbb{T}_r$ from the set $\{k+2, \dots, n+k+1\}$
to be partitioned into new blocks, and distribute the remainder $n-r$ elements among
the $j-1$ blocks $\mathbb{S}_2,\dots,\mathbb{S}_j$, then merge all the singletons formed
from the $r$-subset $\mathbb{T}_r$ (having been partitioned) into $\mathbb{S}_1$ to form one block.
There are $\binom{n}{r}$ ways to choose the $r$ elements, $B_r$
ways to partition them into new blocks, and $(j-1)^{n-r}$ ways to distribute the remainder $n-r$
elements among the $j-1$ blocks. Thus there are $\binom{n}{r}B_{r}S(k,j)(j-1)^{n-r}$ of such partitions. Note that $k+1$ is the largest
singleton, summing over all possible values of $j$ and $r$ produces all ways to partition the set $[n+k+1]$ with the
largest singleton $k+1$. This gives a proof of (\ref{eqn 2.7}). \qed\vskip.2cm

\section{Identities involving $A_{n,k}$ and Bell numbers $B_n$}

\begin{theorem} For any integer $n\geq 0$ and any indeterminant $y$, there hold
\begin{eqnarray}\label{eqn 3.1}
\sum_{k=0}^{n}(-1)^{n-k}\binom{n}{k}A_{n,k}(y+1)^{k}=\sum_{k=0}^{n}\binom{n}{k}y^{k}B_{k},
\end{eqnarray}
or equivalently
\begin{eqnarray}\label{eqn 3.2}
\sum_{k=0}^{n}\binom{n}{k}A_{n,k}y^{k}=\sum_{k=0}^{n}(-1)^{n-k}\binom{n}{k}(y+1)^{k}B_{k}.
\end{eqnarray}
\end{theorem}
\pf Note that $A(-x,x(y+1))=B(xy)e^{x}$ and $A(x,xy)=B(x(y+1))e^{-x}$ from Lemma \ref{lemma 2.1}, by equating the coefficients of
$\frac{x^n}{n!}$ in the resulting series,
one can easily deduce (\ref{eqn 3.1})-(\ref{eqn 3.2}). Also (\ref{eqn 3.2}) can be obtained from (\ref{eqn 3.1}) by setting $y:=-y-1$.
One can be asked to give a combinatorial proof for these two identities. \qed\vskip.2cm

\begin{corollary} For any integer $n\geq 0$, there hold
\begin{eqnarray}
\sum_{k=0}^{n}(-1)^{n-k}\binom{n}{k}A_{n,k}      &=& 1 ,  \label{eqn 3.3}\\
\sum_{k=0}^{n}(-1)^{n-k}\binom{n}{k}2^{k}A_{n,k} &=& B_{n+1}. \label{eqn 3.4}
\end{eqnarray}
\end{corollary}
\pf The case $y=0$ in (\ref{eqn 3.1}) yields (\ref{eqn 3.3}). The case $y=1$ in (\ref{eqn 3.1}), together with (\ref{eqn 1.2}),
yields (\ref{eqn 3.4}). \qed\vskip0.2cm

\begin{corollary} For any integer $n\geq 0$ and any indeterminant $y$, there hold
\begin{eqnarray}
\sum_{k=0}^{n}(-1)^{n-k}\binom{n}{k}A_{n,k}B_{k+1}(y) &=& y\sum_{k=0}^{n}\binom{n}{k}B_{k}B_k(y),  \label{eqn 3.5}\\
\sum_{k=0}^{n}(-1)^{n-k}\binom{n}{k}B_{k}B_{k+1}(y) &=& y\sum_{k=0}^{n}\binom{n}{k}A_{n,k}B_{k}(y). \label{eqn 3.6}
\end{eqnarray}
\end{corollary}
\pf This is an equivalent form of Theorem 3.1. Define a linear (invertible) transformation
\begin{eqnarray*}
L_1(y^k)=B_k(y),\ \ \ (k=0,1,2,\dots).
\end{eqnarray*}
It is well known that $B_k(y)$ satisfies the relation
\begin{eqnarray*}
B_{n+1}(y)=y\sum_{k=0}^{n}\binom{n}{k}B_k(y).
\end{eqnarray*}
Then we have
\begin{eqnarray*}
yL_1((y+1)^n)=y\sum_{k=0}^{n}\binom{n}{k}L_1(y^k)=y\sum_{k=0}^{n}\binom{n}{k}B_k(y)=B_{n+1}(y).
\end{eqnarray*}
Hence (\ref{eqn 3.5}) and (\ref{eqn 3.6}) follow by acting $yL_1$ on the two sides of (\ref{eqn 3.1}) and (\ref{eqn 3.2})
respectively.  \qed\vskip0.2cm
Similarly, if define another linear transformation
\begin{eqnarray*}
L_2(y^k)=\binom{y}{k},\ \ \ (k=0,1,2,\dots),
\end{eqnarray*}
by the Vandermonde's convolution identity
\begin{eqnarray*}
\sum_{k=0}^{n}\binom{a}{k}\binom{b}{n-k}=\binom{a+b}{n},
\end{eqnarray*}
we have
\begin{eqnarray*}
L_2((y+1)^n)=\sum_{k=0}^{n}\binom{n}{k}L_2(y^k)=\binom{y+n}{n}.
\end{eqnarray*}
Then acting $L_2$ on the two sides of (\ref{eqn 3.1}) and (\ref{eqn 3.2}) leads respectively to another equivalent form of Theorem 3.1.
\begin{corollary} For any integer $n\geq 0$ and any indeterminant $y$, there hold
\begin{eqnarray*}
\sum_{k=0}^{n}(-1)^{n-k}\binom{n}{k}\binom{y+k}{k}A_{n,k} &=& \sum_{k=0}^{n}\binom{n}{k}\binom{y}{k}B_{k}, \\
\sum_{k=0}^{n}(-1)^{n-k}\binom{n}{k}\binom{y+k}{k}B_{k}   &=& \sum_{k=0}^{n}\binom{n}{k}\binom{y}{k}A_{n,k}.
\end{eqnarray*}
\end{corollary}
With the Bell umbra $\mathbf{B}$ \cite{Gesselb, Roman, RomRota}, given by $\mathbf{B}^{n}=B_{n}$, (\ref{eqn 1.2}) may be written as
$\mathbf{B}^{n+1}=(\mathbf{B}+1)^{n}$. By (\ref{eqn 2.4}), $A_{n,k}$ can be written umbrally as
\begin{eqnarray*}
A_{n,k}=\mathbf{B}^{k}(\mathbf{B}-1)^{n-k}.
\end{eqnarray*}
Setting $y=\frac{y}{1-y}$ in (\ref{eqn 3.1}) and (\ref{eqn 3.2}), and multiplying $y^m(1-y)^n$ by their two sides, we have
\begin{eqnarray*}
\sum_{k=0}^{n}\binom{n}{k}A_{n,k}y^m(y-1)^{(n+m-k)-m} &=& \sum_{k=0}^{n}(-1)^{n-k}\binom{n}{k}y^{m+k}(y-1)^{(n+m)-(m+k)}B_{k}, \\
\sum_{k=0}^{n}\binom{n}{k}y^m(y-1)^{(n+m-k)-m}B_{k} &=& \sum_{k=0}^{n}(-1)^{n-k}\binom{n}{k}y^{m+k}(y-1)^{(n+m)-(m+k)}A_{n,k},
\end{eqnarray*}
which, when $y=\mathbf{B}$, produce another two identities.
\begin{corollary} For any integers $n,m\geq 0$, there hold
\begin{eqnarray*}
\sum_{k=0}^{n}(-1)^{n-k}\binom{n}{k}A_{n+m,m+k}B_{k}    &=& \sum_{k=0}^{n}\binom{n}{k}A_{n+m-k,m}A_{n,k}, \\
\sum_{k=0}^{n}(-1)^{n-k}\binom{n}{k}A_{n+m,m+k}A_{n,k}  &=& \sum_{k=0}^{n}\binom{n}{k}A_{n+m-k,m}B_{k}.
\end{eqnarray*}
\end{corollary}

\begin{theorem} For any integers $n,k\geq 0$ and any indeterminant $y$, there hold
\begin{eqnarray}
\sum_{j=0}^{n}\binom{n}{j}A_{k+j,k}(y+1)^{n-j}  &=& \sum_{j=0}^{n}\binom{n}{j}B_{k+j}y^{n-j}, \label{eqn 3.7}\\
\sum_{j=0}^{n}\binom{n}{j}A_{k+j,k}B_{n-j+1}(y)  &=& y\sum_{j=0}^{n}\binom{n}{j}B_{k+j}B_{n-j}(y), \label{eqn 3.8}\\
\sum_{j=0}^{n}\binom{n}{j}\binom{y+n-j}{n-j}A_{k+j,k}  &=& \sum_{j=0}^{n}\binom{n}{j}\binom{y}{n-j}B_{k+j}. \label{eqn 3.9}
\end{eqnarray}
\end{theorem}
\pf Note that $A(x,t)e^{x(y+1)}=B(x+t)e^{xy}$ from Lemma \ref{lemma 2.1}, by equating the coefficients of
$\frac{x^nt^k}{n!k!}$ in the resulting series,
one can easily deduce (\ref{eqn 3.7}). (\ref{eqn 3.8}) and (\ref{eqn 3.9}) can be followed respectively by acting $yL_1$ and
$L_2$ on the two sides of (\ref{eqn 3.7}). Here we provide a combinatorial proof for (\ref{eqn 3.7}).

Let $\mathbb{X}_{n,k}=\bigcup_{j=0}^n\mathbb{X}_{n, k, j}$ and $\mathbb{X}_{n, k, j}$ denote the set of pairs $(\pi, \mathbb{S})$ such that
\begin{itemize}
\item $\mathbb{S}$ is an $(n-j)$-subset of $[k+2, n+k+1]=\{k+2, \dots, n+k+1\}$, and each element of $\mathbb{S}$ has weight $1$ or $y$; In
other words, each element of $\mathbb{S}$ has weight $1+y$;
\item $\pi$ is a partition of the set $[n+k+1]-\mathbb{S}$ with the largest singleton $k+1$, and each element of $[n+k+1]-\mathbb{S}$ has weight $1$.
\end{itemize}

Let $\mathbb{Y}_{n,k}=\bigcup_{j=0}^n\mathbb{Y}_{n, k, j}$ and $\mathbb{Y}_{n, k, j}$ denote the set of pairs $(\pi, \mathbb{S})$ such that
\begin{itemize}
\item $\mathbb{S}$ is an $(n-j)$-subset of $[k+2, n+k+1]$ and each element of $\mathbb{S}$ has weight $y$;
\item $\pi$ is a partition of the set $[n+k+1]-\mathbb{S}$ such that $k+1$ must be a singleton, and each element of $[n+k+1]-\mathbb{S}$ has weight $1$.
\end{itemize}
The weight of $(\pi, \mathbb{S})$ is defined to be the product of the weight of each element of $[n+k+1]$. Clearly, the weights of $\mathbb{X}_{n,k}$
and $\mathbb{Y}_{n,k}$ are counted respectively by the left and right sides of (\ref{eqn 3.7}).

Given any pair $(\pi, \mathbb{S})\in \mathbb{X}_{n,k}$, $\mathbb{S}$ can be partitioned into two parts $\mathbb{S}_1$ and $\mathbb{S}_{2}$ such that each
element of $\mathbb{S}_1$ has weight $y$ and each element of $\mathbb{S}_2$ has weight $1$. Regard each element of $\mathbb{S}_2$ as a singleton,
together with $\pi$, we obtain a partition $\pi_1$ of $[n+k+1]-\mathbb{S}_1$ such that $k+1$ is a singleton.
Then the pair $(\pi_1, \mathbb{S}_1)$ lies in $\mathbb{Y}_{n,k}$.

Conversely, for any pair $(\pi_1, \mathbb{S}_1)\in \mathbb{Y}_{n,k}$, let $\mathbb{S}$ denote the union of $\mathbb{S}_1$ and the singletons of $\pi_1$
greater than $k+1$, then $\pi_1$ can be partitioned into two parts $\pi$ and $\pi'$ such that $\pi$ is a partition of $[n+k+1]-\mathbb{S}$ with
the largest singleton $k+1$ and $\pi'$ is the singletons of $\pi_1$ greater than $k+1$. Then the pair $(\pi, \mathbb{S})$ lies in $\mathbb{X}_{n,k}$.

Clearly we find a bijection between $\mathbb{X}_{n,k}$ and $\mathbb{Y}_{n,k}$, which proves (\ref{eqn 3.7}). \qed\vskip.2cm

Setting $y=0$ and $y=1$ in (\ref{eqn 3.7}), by (\ref{eqn 2.5}) in the case $k=1$, we have
\begin{corollary} For any integers $n,k \geq 0$, there hold
\begin{eqnarray*}
B_{n+k}       &=& \sum_{j=0}^{n}\binom{n}{j}A_{k+j,k},  \\
A_{n+k+1,n+1} &=& \sum_{j=0}^{n}\binom{n}{j}A_{k+j,k}2^{n-j}.
\end{eqnarray*}
\end{corollary}

\begin{corollary} For any integers $n,k,m,i\geq 0$, there hold
\begin{eqnarray}
\sum_{j=0}^{n}\binom{n}{j}A_{k+j,k}(n-j)!  &=& \sum_{j=0}^{n}\binom{n}{j}B_{k+j}D_{n-j},  \label{eqn 3.10} \\
\sum_{j=0}^{n}(-1)^{n-j}\binom{n}{j}A_{k+j,k}A_{m+i+j,m}  &=& \sum_{j=0}^{n}(-1)^{n-j}\binom{n}{j}A_{n+m+i,n+m-j}B_{k+j}, \label{eqn 3.11}
\end{eqnarray}
where $D_n$ is the number of permutations of $[n]$ without fixed points.
\end{corollary}
\pf The exponential generating function \cite{Stanley} for $D_n$ is
\begin{eqnarray*}
\sum_{n\geq 0}D_n\frac{x^n}{n!}=\frac{e^{-x}}{1-x},
\end{eqnarray*}
from which, one can get
\begin{eqnarray*}
n!=\sum_{j=0}^{n}\binom{n}{j}D_{n-j}.
\end{eqnarray*}
Let $\mathbf{D}$ be the umbra, given by $\mathbf{D}^n=D_n$, we have $n!=(\mathbf{D}+1)^n$. Then (\ref{eqn 3.10}) can be obtained by setting
$y=\mathbf{D}$ in (\ref{eqn 3.7}).

Setting $y=\frac{y}{1-y}$ in (\ref{eqn 3.7}) and multiplying $y^m(y-1)^{n+i}$ by the two sides, we have
\begin{eqnarray*}
\sum_{j=0}^{n}(-1)^{n-j}\binom{n}{j}A_{k+j,k}y^m(y-1)^{i+j}  &=& \sum_{j=0}^{n}(-1)^{n-j}\binom{n}{j}B_{k+j}y^{n+m-j}(y-1)^{i+j},
\end{eqnarray*}
which, when $y=\mathbf{B}$, yields (\ref{eqn 3.11}). \qed\vskip.2cm

Gould and Quaintance \cite{GouldQuain} present the identity
\begin{eqnarray*}
\sum_{j=0}^{m}s(m,j)B_{k+j} &=& \sum_{i=0}^{k}\binom{k}{i}m^{k-i}B_{i},
\end{eqnarray*}
which is a special case ($n=0$) of the following three identities.

\begin{theorem} For any integers $n,m,k\geq 0$, there hold
\begin{eqnarray}
\sum_{j=0}^{m}s(m,j)A_{n+k+j,k+j} &=& \sum_{i=0}^{k}\sum_{r=0}^{n}\binom{k}{i}\binom{n}{r}m^{n+k-i-r}A_{r+i,i}, \label{eqn 3.12} \\
\sum_{j=0}^{m}s(m,j)A_{n+k+j,k+j} &=& \sum_{i=0}^{k}\sum_{r=0}^{n}\binom{k}{i}\binom{n}{r}m^{k-i}(m-1)^{n-r}B_{r+i}, \label{eqn 3.13} \\
\sum_{j=0}^{m}s(m,j)A_{n+k+j+1,n+1} &=& \sum_{i=0}^{k}\sum_{r=0}^{n}\binom{k}{i}\binom{n}{r}m^{k-i}(m+1)^{n-r}B_{r+i}, \label{eqn 3.14}
\end{eqnarray}
where $s(k,j)$ are the first kind of Stirling numbers.
\end{theorem}
\pf We know the Bell umbra $\mathbf{B}$ satisfies $\mathbf{B}^{m+1}=(\mathbf{B}+1)^{m}$. Then
by linearity, for any polynomial $f(x)$ we have
\begin{eqnarray*}
\mathbf{B}f(\mathbf{B})=f(\mathbf{B}+1),
\end{eqnarray*}
which, by induction on integer $m\geq 0$, leads to
\begin{eqnarray}\label{eqn 3.15}
\mathbf{B}(\mathbf{B}-1)\cdots(\mathbf{B}-m+1)f(\mathbf{B})=f(\mathbf{B}+m).
\end{eqnarray}

It is well known that for any indeterminant $x$,
\begin{eqnarray*}
x(x-1)\cdots (x-m+1)=\sum_{j=0}^{m}s(m,j)x^j,
\end{eqnarray*}
Using the umbral representation for $A_{n,k}$, we have
\begin{eqnarray*}
\sum_{j=0}^{m}s(m,j)A_{n+k+j,k+j}&=& \sum_{j=0}^{m}s(m,j)\mathbf{B}^{k+j}(\mathbf{B}-1)^{n}       \\
                                 &=& \mathbf{B}(\mathbf{B}-1)\cdots(\mathbf{B}-m+1)\mathbf{B}^{k}(\mathbf{B}-1)^{n}      \\
                                 &=& (\mathbf{B}+m)^{k}(\mathbf{B}-1+m)^{n}            \\
                                 &=& \sum_{i=0}^{k}\sum_{r=0}^{n}\binom{k}{i}\binom{n}{r}m^{n+k-i-r}\mathbf{B}^{i}(\mathbf{B}-1)^{r} \\
                                 &=& \sum_{i=0}^{k}\sum_{r=0}^{n}\binom{k}{i}\binom{n}{r}m^{n+k-i-r}A_{r+i,i},
\end{eqnarray*}
which proves (\ref{eqn 3.12}). Similarly, one can deduce (\ref{eqn 3.13}) and (\ref{eqn 3.14}). \qed\vskip.2cm

\begin{theorem} For any integer $n\geq 0$, there hold
\begin{eqnarray}
\sum_{k=0}^{n}(k+1)A_{n,k}  &=& (n+2)B_{n+1}-V_{n+3}, \label{eqn 3.16} \\
\sum_{k=0}^{n}(n-k+1)A_{n,k}&=& V_{n+3}-(n+2)V_{n+1} . \label{eqn 3.17}
\end{eqnarray}
\end{theorem}
\pf Define
\begin{eqnarray*}
\alpha_n(x)=\sum_{k=0}^{n}A_{n,k}x^{k}.
\end{eqnarray*}
By (\ref{eqn 2.2}), we have
\begin{eqnarray*}
 \sum_{k=0}^{n}A_{n,k}x^{k}&=& \frac{1}{e}\sum_{k=0}^{n}x^k\sum_{m=0}^{\infty}\frac{m^k(m-1)^{n-k}}{m!}     \\
                           &=& \frac{1}{e}\sum_{m=0}^{\infty}\frac{1}{m!}\sum_{k=0}^{n}(mx)^k(m-1)^{n-k}     \\
                           &=& \frac{1}{e}\sum_{m=0}^{\infty}\frac{1}{m!}\frac{(mx)^{n+1}-(m-1)^{n+1}}{mx-m+1}.  \\
\end{eqnarray*}
Differentiating $x\alpha_n(x)$ and then setting $x=1$ gives
\begin{eqnarray*}
\sum_{k=0}^{n}(k+1)A_{n,k} &=& \frac{1}{e}\sum_{m=0}^{\infty}\frac{(n+1)m^{n+1}-m^{n+1}(m-1)+(m-1)^{n+2}}{m!} \\
                           &=& (n+1)B_{n+1}-A_{n+2,n+1}+A_{n+2,0} \\
                           &=& (n+1)B_{n+1}-(A_{n+2,n+2}-A_{n+1,n+1})+(B_{n+2}-A_{n+3,0}) \\
                           &=& (n+1)B_{n+1}-(B_{n+2}-B_{n+1})+(B_{n+2}-V_{n+3}) \\
                           &=& (n+2)B_{n+1}-V_{n+3},
\end{eqnarray*}
which proves (\ref{eqn 3.16}). Similarly, differentiating $x^{n+1}\alpha_n(x^{-1})$ and then setting $x=1$ gives (\ref{eqn 3.17}). \qed\vskip.2cm

\begin{remark}
Canfield \cite{Canfield} has shown that the average number of singletons in a partition of $[n]$ is
an increasing function of $n$. We guess that the average number of the largest or smallest singletons in a partition of $[n+1]$
is also an increasing function of $n$. That is to say, both
\begin{eqnarray*}
\frac{(n+2)B_{n+1}-V_{n+3}}{B_{n+1}}\hskip.5cm {\mbox{and}} \hskip.5cm \frac{V_{n+3}-(n+2)V_{n+1}}{B_{n+1}}
\end{eqnarray*}
are increasing functions of $n$. One can be asked for asymptotic formulas for the above two expressions.
\end{remark}

\section{Congruence properties of $A_{n,k}$ and Bell numbers $B_n$}

In this section, based on umbral calculus, we study the congruence properties of $A_{n,k}$ and Bell numbers $B_n$. Throughout this section,
$p$ refers to a prime, and unless stated otherwise, all congruences are modulo $p$.

\begin{theorem} For any integers $n, m, k\geq 0$, there holds
\begin{eqnarray*}
A_{n+pm+k, k}\equiv A_{n+m+k, m+k}.
\end{eqnarray*}
\end{theorem}
\pf Recall that the Lagrange congruence
\begin{eqnarray*}
x(x-1)\cdots (x-p+1) &\equiv& x^{p}-x,
\end{eqnarray*}
and the binomial congruence
\begin{eqnarray*}
(x-1)^{p} &\equiv& x^{p}-1.
\end{eqnarray*}
Setting $x=\mathbf{B}$, by (\ref{eqn 3.15}), for any polynomial $f(x)$, one gets
\begin{eqnarray*}
(\mathbf{B}^p-\mathbf{B})f(\mathbf{B}) \equiv \mathbf{B}(\mathbf{B}-1)\cdots(\mathbf{B}-p+1)f(\mathbf{B})
=f(\mathbf{B}+p)\equiv f(\mathbf{B}),
\end{eqnarray*}
which, by induction on integer $j\geq 0$, leads to
\begin{eqnarray}
(\mathbf{B}^p-\mathbf{B})^{j}f(\mathbf{B}) &\equiv & f(\mathbf{B}),  \label{congru A}  \\
(\mathbf{B}-1)^{pj}f(\mathbf{B}) &\equiv & \mathbf{B}^{j}f(\mathbf{B}). \label{congru B}
\end{eqnarray}

Using the umbral representation for $A_{n,k}$, we have
\begin{eqnarray*}
A_{n+pm+k,k} &=&       \mathbf{B}^{k}(\mathbf{B}-1)^{n+pm}     \\
             &=&       (\mathbf{B}-1)^{pm}\mathbf{B}^{k}(\mathbf{B}-1)^{n}    \\
             &\equiv & \mathbf{B}^{m+k}(\mathbf{B}-1)^{n}       \\
             &=&       A_{n+m+k,m+k},
\end{eqnarray*}
as claimed. \qed\vskip.2cm

\begin{corollary} For any integers $n, m\geq 0$, there hold
\begin{eqnarray}
B_{n+pm}      & \equiv & A_{n+m+1, m+1},  \label{eqn 4.1} \\
B_{n+p}       & \equiv & B_n+B_{n+1},  \hskip1cm  \mbox{(Touchard's\ congruence\ \cite{Toucha, Touchb})}, \label{eqn 4.2} \\
A_{(n+1)p, p} & \equiv & B_{n}+B_{n+1}, \label{eqn 4.3}  \\
B_{np}        & \equiv & B_{n+1}, \hskip2cm  {\mbox{(Comtet's\ congruence\ \cite{Comtet, GertRob})}}. \label{eqn 4.4}
\end{eqnarray}
\end{corollary}
\pf The case $k=1$ in Theorem 4.1 leads to (\ref{eqn 4.1}), which in the case $m=1$ yields (\ref{eqn 4.2}). (\ref{eqn 4.3})
follows by setting $n=0, m=n, k=p$, and (\ref{eqn 4.4}) follows by setting $n=0, m=n, k=1$ in Theorem 4.1. \qed\vskip.2cm

\begin{theorem} For any integers $n, m, k\geq 0$, there holds
\begin{eqnarray*}
A_{n+p^m+k, k} &\equiv & mA_{n+k,k}+A_{n+k+1,k}.
\end{eqnarray*}
\end{theorem}
\pf By (\ref{congru A}) and (\ref{congru B}), when $f(x)=1$, one has
\begin{eqnarray*}
\mathbf{B}^{p} & \equiv & \mathbf{B}+1, \\
(\mathbf{B}-1)^{p} & \equiv & \mathbf{B}.
\end{eqnarray*}
Using the little Fermat's congruence $k^p  \equiv  k$, where $k$ is an integer,
by induction on integer $m\geq 0$, we have
\begin{eqnarray}\label{congru C}
(\mathbf{B}-1)^{p^m} & \equiv & \mathbf{B}+m-1.
\end{eqnarray}
Then
\begin{eqnarray*}
A_{n+p^m+k, k} &= &          (\mathbf{B}-1)^{p^m}\mathbf{B}^{k}(\mathbf{B}-1)^{n}  \\
               &\equiv &     (\mathbf{B}+m-1)\mathbf{B}^{k}(\mathbf{B}-1)^{n}  \\
               &= &           m\mathbf{B}^{k}(\mathbf{B}-1)^{n}+\mathbf{B}^{k}(\mathbf{B}-1)^{n+1} \\
               &=&            mA_{n+k,k}+A_{n+k+1,k},
\end{eqnarray*}
as desired. \qed\vskip.2cm

\begin{theorem} Let $N_p=\frac{p^p-1}{p-1}$, for any integers $n, k\geq 0$, there hold
\begin{eqnarray*}
A_{n+N_p+k, k} &\equiv & A_{n+k,k}, \\
A_{n+N_p+k, N_p+k} &\equiv & A_{n+k,k},
\end{eqnarray*}
namely, the sequences $(A_{n+k,k})_{n\geq 0}$ and $(A_{n+k,k})_{k\geq 0}$ (mod $p$) both have the period $N_p$.
\end{theorem}
\pf By (\ref{congru C}) and the Lagrange cogruence, one has
\begin{eqnarray*}
(\mathbf{B}-1)^{N_p}=\prod_{j=1}^{p}(\mathbf{B}-1)^{p^{p-j}}  \equiv   \prod_{j=1}^{p}(\mathbf{B}-j-1)
\equiv   \prod_{j=0}^{p-1}(\mathbf{B}-j)\equiv 1.
\end{eqnarray*}
Then
$$A_{n+N_p+k, k}= (\mathbf{B}-1)^{N_p}\mathbf{B}^{k}(\mathbf{B}-1)^{n}
\equiv \mathbf{B}^{k}(\mathbf{B}-1)^{n} =A_{n+k,k}.$$

When $m=N_p$ in Theorem 4.1, one has
$$A_{n+N_p+k, N_p+k}\equiv A_{n+pN_p+k, k}\equiv  A_{n+k,k},$$
where the last modular equation follows by the periodicity of $(A_{n+k,k})_{n\geq 0}$. \qed\vskip.2cm

\begin{remark}
Hall showed that the Bell numbers (the case $k=1$ for $(A_{n+k,k})_{n\geq 0}$ or the case $n=0$
for $(A_{n+k,k})_{k\geq 0}$) have the period $N_p$, a result rediscovered by
Williams \cite{Williams}. Williams also showed that the minimum period is precisely $N_p$ for $p = 2, 3$ and $5$.
Radoux \cite{Radoux} conjectured that $N_p$ is the minimal period of the sequence $B_n$ for any prime $p$.
Levine and Dalton \cite{LevDal} showed that the minimum period is exactly $N_p$ for $p = 7, 11, 13$
and $17$. They also investigated the period for the other primes $< 50$. Recently, Montgomery, Nahm and Wagstaff \cite{MonNahWag} showed
that the minimum period is exactly $N_p$ for most primes $p$ below 180.
\end{remark}

For the sequences$(A_{n+k,k})_{n\geq 0}$ and $(A_{n+k,k})_{k\geq 0}$, we also have the following conjecture.

\begin{conjecture} For any integer $k\geq 0$ and any prime $p$, the sequences $(A_{n+k,k})_{n\geq 0}$ and $(A_{n+k,k})_{k\geq 0}$ both have
the minimum period $N_p$ modulo $p$.
\end{conjecture}

\begin{theorem} Let $n,m, k\geq 0$ be integers and $p$ be a prime. Then a necessary and sufficient condition that
$A_{n+m+k,k}\equiv 0\ (\mbox{mod}\ p)$ for $m=0, 1, \dots, p-2$, is that $A_{n+m+k,k}\equiv A_{n+pm+k,k} \ (\mbox{mod}\ p)$
for $m=1, 2, \dots, p-1$.
\end{theorem}
\pf By Theorem 4.1 and (\ref{eqn 2.5}), we have
\begin{eqnarray}\label{eqn 4.5}
A_{n+pm+k, k} &\equiv & A_{n+m+k,m+k}=\sum_{j=0}^m\binom{m}{j}A_{n+k+j,k}.
\end{eqnarray}
Therefore, if $A_{n+m+k,k}\equiv 0\ (\mbox{mod}\ p)$ for $m=0, 1, \dots, p-2$, we clearly have $A_{n+pm+k, k} \equiv 0$ and
hence, trivially, $A_{n+pm+k, k}  \equiv  A_{n+m+k,k}(\equiv 0)$. When $m=p-1$ and $A_{n+j+k,k}\equiv 0 $ for $j=0, 1, \dots, p-2$,
(\ref{eqn 4.5}) reduces to $A_{n+(p-1)+k,k}\equiv A_{n+p(p-1)+k,k}$.

Conversely, if $A_{n+m+k,k}\equiv A_{n+pm+k,k} \ (\mbox{mod}\ p)$ for $m=1, 2, \dots, p-1$, (\ref{eqn 4.5}) is equivalent to
\begin{eqnarray*}
A_{n+m+k, k} &\equiv & A_{n+m+k,k}+\sum_{j=0}^{m-1}\binom{m}{j}A_{n+k+j,k}, \ (m=1,2,\dots, p-1),
\end{eqnarray*}
which reduces to
\begin{eqnarray}\label{eqn 4.6}
0 &\equiv & \sum_{j=0}^{m-1}\binom{m}{j}A_{n+k+j,k}, \ (m=1,2,\dots, p-1).
\end{eqnarray}

The system (\ref{eqn 4.6}) is triangular with diagonal coefficients $\binom{m}{m-1}$. The coefficient matrix is therefore nonsingular with
determinant $(p-1)!\equiv -1$ by Wilson's theorem. Thus the only solution is given by
$A_{n+m+k,k}\equiv 0\ (\mbox{mod}\ p)$ for $m=0, 1, \dots, p-2.$ \qed\vskip.2cm

\begin{theorem} For any integer $k\geq 0$ and any prime $p$, there exists an integer $M_{p,k}\geq 0$ such that
\begin{eqnarray*}
A_{M_{p,k}+m+k, k} &\equiv & 0, \ (0\leq m\leq p-2),
\end{eqnarray*}
where
\begin{eqnarray*}
 M_{p,k}  &\equiv & 1-(k-1)p-\frac{p^p-p}{(p-1)^2}, \hskip1cm (\mbox{mod}\ N_p).
\end{eqnarray*}
In other words, the sequence $(A_{n+k,k})_{n\geq 0} \ (\mbox{mod} \ p)$ contains a string of $p-1$ consecutive zeroes.
\end{theorem}
\pf By Theorem 4.1 and (\ref{eqn 4.1}), we have
\begin{eqnarray*}
A_{(n+(k-1)p-k)p+k, k} &\equiv & A_{n+k,k},
\end{eqnarray*}
which, when $n=M_{p,k}+m$, where $M_{p,k}$ is an integer to be determined, produces
\begin{eqnarray*}
A_{(M_{p,k}+(k-1)p-k)p+pm+k, k} &\equiv & A_{M_{p,k}+m+k,k}.
\end{eqnarray*}

By Theorem 4.4 and 4.7, it follows that $p-1$ consecutive zeros of $(A_{n+k,k})_{n\geq 0} \ (\mbox{mod} \ p)$ will occur,
beginning with $A_{M_{p,k}+k,k}$, if there holds
\begin{eqnarray*}
 A_{(M_{p,k}+(k-1)p-k)p+pm+k,k}  &\equiv & A_{(M_{p,k}+(k-1)p-k)p+m+k,k}, \ (m=1,2,\dots, p-1).
\end{eqnarray*}
It is just required that the following condition holds
\begin{eqnarray*}
 (M_{p,k}+(k-1)p-k)p+m  &\equiv & M_{p,k}+m,  \hskip1cm (\mbox{mod}\ N_p),
\end{eqnarray*}
or, equivalently, if there holds
\begin{eqnarray*}
 (M_{p,k}+(k-1)p-k)p   &=& M_{p,k}+rN_p,
\end{eqnarray*}
for some integer $r$. Using $N_p=\frac{p^p-1}{p-1}=\frac{p^p-p}{p-1}+1$, routine calculation yields
\begin{eqnarray*}
M_{p,k}  &=& 1-(k-1)p+\frac{r+1}{p-1}+r\frac{p^p-p}{(p-1)^2},
\end{eqnarray*}
It is easy to verify by the binomial congruence that $\frac{p^p-p}{(p-1)^2}$ is always an integer. Since $M_{p,k}$ is also
an integer, so we must have $r=-1+t(p-1)$ for some integer $t$, from which it follows that
\begin{eqnarray*}
M_{p,k}  &=& 1-(k-1)p-\frac{p^p-p}{(p-1)^2}+tN_p.
\end{eqnarray*}
Since the $p-1$ consecutive zeros start with $A_{M_{p,k}+k,k}$, the proof is complete. \qed\vskip.2cm

Using the same arguments, we have analogous results for the sequences $(A_{n+k,n})_{n\geq 0}$, their proofs are left to interested readers,
the critical step for Theorem 4.10 is to show the congruence relation
\begin{eqnarray*}
A_{(n-1)p-k(p^{p-1}-1)+k, (n-1)p-k(p^{p-1}-1)} &\equiv & A_{n+k,n}.
\end{eqnarray*}

\begin{theorem} Let $n,m, k\geq 0$ be integers and $p$ be a prime. Then a necessary and sufficient condition that
$A_{n+m+k,n+m}\equiv 0\ (\mbox{mod}\ p)$ for $m=0, 1, \dots, p-2$, is that $A_{n+m+k,n+m}\equiv A_{n+pm+k,n+pm} \ (\mbox{mod}\ p)$
for $m=1, 2, \dots, p-1$.
\end{theorem}

\begin{theorem} For any integer $k\geq 0$ and any prime $p$, there exists an integer $U_{p,k}\geq 0$ such that
\begin{eqnarray*}
A_{U_{p,k}+m+k, U_{p,k}+m} &\equiv & 0, \ (0\leq m\leq p-2),
\end{eqnarray*}
where
\begin{eqnarray*}
 U_{p,k}  &\equiv & 1+\frac{(p^{p-1}-1)k}{p-1}-\frac{p^p-p}{(p-1)^2}, \hskip1cm (\mbox{mod}\ N_p).
\end{eqnarray*}
In other words, the sequence $(A_{n+k,n})_{n\geq 0} \ (\mbox{mod} \ p)$ contains a string of $p-1$ consecutive zeroes.
\end{theorem}

\begin{remark}
Radoux \cite{Radoux} shows that if the period of the residues of the Bell sequence $B_n$ is equal to $N_p$ for a given prime
$p$, then there exists a number of $c$, depending on $p$, such that $B_{c+m}\equiv 0$ (mod $p$) for $0\leq m\leq p-2$. He also obtains the location
of such a string of consecutive zeros. Kahale \cite{Kahale} and Layman \cite{Layman} show respectively by two entirely different methods that this result
holds without the hypothesis that $N_p$ is the minimal period. Their result is a special case of Theorem 4.8 for $k=1$ or of Theorem 4.10 for $k=0$.
Our proof methods are similar to Layman's.
\end{remark}

\section*{Acknowledgements} The authors are grateful to Eva YuPing Deng and the
anonymous referees for the helpful suggestions and comments. The
work was supported by The National Science Foundation of China
(Grant No. 10801020 and 70971014).


\end{document}